\documentclass[12pt,twoside]{article}
\usepackage{amsmath,amsbsy,amsfonts}
\usepackage{graphicx,color}
\def\dis{\displaystyle}
\def\nd{\noindent}
\def\thend{\rule{3mm}{3mm}}
\def\Re{\mathbb{R}}

\newtheorem{thm}{Theorem}[section]

\newtheorem{lem}{Lemma}[section]
\newtheorem{rmk}{Remark}[section]

\newtheorem{definition}{Definition}[section]
\newcommand{\fim}{\hfill\rule{2mm}{2mm}}
\newenvironment{dem}{\nd{\bf Proof. }}{\hskip.3cm\thend}

\begin{document}

\setlength{\baselineskip}{6.5mm} \setlength{\oddsidemargin}{8mm}
\setlength{\topmargin}{-3mm}
\title{\Large\sf A Sub-Supersolution Approach for a \\ Quasilinear Kirchhoff Equation}
\author{\sf
 Claudianor O. Alves\thanks{Partially supported by CNPq - Grant 304036/2013-7} \;\; and \;\;  Francisco Julio S.A. Corr\^ea \thanks{Partially supported by CNPq - Grant 301807/2013-2 } \\
 Universidade Federal de Campina Grande\\
 Unidade Acad\^emica de Matem\'atica - UAMat\\
 58.429-900 - Campina Grande - PB - Brazil\\
 coalves@dme.ufcg.edu.br \& fjsacorrea@pq.cnpq.br  }

\pretolerance10000
\date{}
\numberwithin{equation}{section} \maketitle

\begin{abstract}
In this paper we establish an  existence result for a quasilinear Kirchhoff equation via a sub and supersolution approach, by using the pseudomonotone operators theory.

\vspace{0.2cm} \noindent \emph{2000 Mathematics Subject
Classification} : 34B18, 35A15, 46E39.


\end{abstract}

\section{Introduction}

In this paper we deal with the quasilinear stationary Kirchhoff equation
$$
\left\{
\begin{array}{rclcc}
  -M\left(\dis\int_{\Omega}|\nabla u|^{2}dx\right)\Delta u & = & f(x,u, \nabla u ) & \mbox{in} & \Omega , \\

  u  & = & 0 & \mbox{on} & \partial \Omega
\end{array}
\right.
\eqno{(P)}
$$
where $\Omega \subset \Re^{N}, N\geq 1$, is a bounded smooth domain, $f:\Omega \times \Re \times \Re^{N}\rightarrow [0,+\infty)$ is a continuous function satisfying

\vspace{0.5cm}

\noindent $(f_{1})$ There are a continuous function $h:\overline{\Omega} \times \mathbb{R} \to [0,+\infty)$ and $\eta  \in  [0,2]$ such that $|f(x,t,y)|\leq h(x,t)(1+|y|^{\eta})$ for all $(x,t,y)\in \Omega \times \mathbb{R} \times \Re^{N}$;

\vspace{0.5cm}
\noindent and $M:[0, +\infty ) \rightarrow [0, +\infty ) $ satisfies

\vspace{0.5cm}

\noindent $(M_{1})$ $M$ is continuous and increasing;

\vspace{0.5cm}

\noindent $(M_{2})$ There is a positive constant $m$ such that $M(t) \geq m>0$ for all $t\in \Re$.

\vspace{0.5cm}

Problem $(P)$ is a generalization of the classical stationary Kirchhoff equation
\begin{equation}\label{kirchhoff classical}
\left\{
\begin{array}{rclcc}
-M\left(\dis\int_{\Omega}|\nabla u|^{2}dx\right)\Delta u & = & f(x,u) & \mbox{in} & \Omega , \\
                                     u & = & 0 & \mbox{on} & \partial\Omega .
\end{array}
\right.
\end{equation}

As it is well known, problem (\ref{kirchhoff classical}) is the general form of the  stationary counterpart of the hyperbolic Kirchhoff equation
\begin{equation}\label{kirchhoff hyperbolic}
\rho \frac{\partial^{2}u}{\partial t^{2}}-\left(\frac{P_{0}}{h}+\frac{E}{2L}\int^{L}_{0}\left|\frac{\partial u}{\partial x}\right|^{2}dx\right)\frac{\partial^{2}u}{\partial x^{2}}=0,
\end{equation}
that appeared at the first time in the work of Kirchhoff \cite{kirchhoff}, in 1883. The equation in (\ref{kirchhoff hyperbolic}) is called Kirchhoff Equation and it extends the classical D'Alembert wave equation, by considering the effects of the changes in the length of the strings during the vibrations.

The interest of the mathematicians on the so called nonlocal problems like (\ref{kirchhoff classical}) (nonlocal because of the presence of the term $M(\|u\|^{2})$, which implies that equations in $(P)$ and (\ref{kirchhoff classical}) are no longer pontwise equalities) has increased because they represent a variety of relevant physical and engineering situations and requires a nontrivial apparatus to solve them. It is worthwhile to emphasize that the  most of the  articles on this subject are concerned with the semilinear case, i.e., $f=f(x,u)$.

In several places we should  face nonhomogeneous Kirchhoff term, that is, the function $M$ also depends on the variable $x\in \Omega$. For instance, L\'imaco, Clark and Medeiros \cite{medeiros} attack a biharmonic evolution equation in which the operator is of the form
$$
\mathcal{L}u\equiv a(x)u''+\Delta (b(x)\Delta u)-M\left( x,t, \int_{\Omega}|\nabla u(x,t)|^{2}dx\right)\Delta u
$$
motivated by the problem of vertical flexion of fully clamped beams. In Figueiredo, Morales-Rodrigo, Santos Junior \& Su\'arez \cite{figueiredo et al} consider a problem whose equation is of the form
$$
-M\left(x, \int_{\Omega}|\nabla u|^{2}dx\right)\Delta u =f(x,u) ~~\mbox{in}~~ \Omega ,
$$
under homogeneous Dirichlet boundary condition, by using a bifurcation argument. Note that for $M$ nonhomogeneous we lose the variational structure and the approach we use in the present article can not be used, at least in a direct way.

In this work, we explore the presence of the gradient term $|\nabla u|$, which makes problem $(P)$ nonvariational, by considering the nonlocal term $M$ with the minimal typical assumptions $(M_{1})-(M_{2})$ which, up to now, at least to our knowledge, has not been considered yet. We point out that in the original Kirchhoff equation the term $M$ is of the form  $M(t)=a+bt, a,b>0$, which  enjoys assumptions $(M_{1})$ and $(M_{2})$.

Our approach was motivated by Cuesta Leon \cite{cuestaleon} and in it the method of sub-supersolution and pseudomonotone operator theory play  a key role. We should say that here we have to surmount several technical difficulties provoked by  the presence of the nonlocal term $M$.

The method of sub and supersolution for semilinear nonlocal equations has been previously used by some authors. We cite some of them.

In Alves-Corr\^ea \cite{alves-correa} the authors study the problem

\begin{equation} \label{ problem1}
\left\{
\begin{array}{rclcc}
  -M\left(\dis\int_{\Omega}|\nabla u|^{2}dx\right)\Delta u & = & f(x,u ) & \mbox{in} & \Omega , \\

  u  & = & 0 & \mbox{on} & \partial \Omega ,
\end{array}
\right.
\end{equation}
via sub-supersolution (monotone iteration) by considering $M:\Re^{+}\rightarrow \Re^{+}$  nonincreasing and $H(t)=M(t^{2})t$  increasing. Note that the typical Kirchhoff term $M(t)=a+bt, a,b>0$ is increasing, i.e., the result in \cite{alves-correa} does not include such a $M$.

In Corr\^ea \cite{correa} the author studies the problem
\begin{equation} \label{correa problem}
\left\{
\begin{array}{rclcc}
  -a\left(\dis\int_{\Omega}|u|^{q}dx\right)\Delta u & = & H(x)f(u ) & \mbox{in} & \Omega , \\

  u  & = & 0 & \mbox{on} & \partial \Omega ,
\end{array}
\right.
\end{equation}
where $a: \Re \rightarrow \Re^{+}$ is a function satisfying $a(s)\geq a_{0}>0 \;  \forall s\in \Re$, $s\mapsto s^{\frac{1}{q}}a(s)$ is increasing and $s\mapsto a(s)$ is decreasing. In particular,  $a$ is a bounded function. In this work the author uses sub-supersolution combined with fixed point theory.

In Chipot-Corr\^ea \cite{chipot-correa} the authors consider the problem

\begin{equation} \label{chipot-correa problem}
\left\{
\begin{array}{rclcc}
  -\mathcal{A}(x,u)\Delta u & = & \lambda f(u ) & \mbox{in} & \Omega , \\

  u  & = & 0 & \mbox{on} & \partial \Omega ,
\end{array}
\right.
\end{equation}
where, among other things, $\mathcal{A}:\Omega \times \Re \rightarrow \Re$ satisfies
\begin{equation}
0<a_{0}\leq \mathcal{A}(x,u) \leq a_{\infty}, \; \mbox{a.e.}\; x\in \Omega , \forall u\in L^{p}(\Omega ).
\end{equation}

In that work, it is used sub-supersolution via fixed point properties and, again, the nonlocal term is bounded.

Here, we permit, inspired by \cite{cuestaleon}, that the Kirchhoff term $M$ may be of the form of the original one.

\begin{definition}\label{weak solution definition}
We say that $u\in H_{0}^{1}(\Omega )\cap L^{\infty}(\Omega )$ is a weak solution of the problem $(P)$ if
\begin{equation}\label{weak solution}
M\left(\int_{\Omega}|\nabla u|^{2}dx\right)\int_{\Omega}\nabla u \nabla vdx=\int_{\Omega}f(x,u, \nabla u)vdx \; \forall v\in H_{0}^{1}(\Omega ).
\end{equation}
\end{definition}

The main result of this paper is as follows:

\begin{thm}\label{main result}
Assume the hypotheses $(M_{1})-(M_{2})$ and $(f_1)$. Moreover, suppose that there are $\overline{u} \in W^{1,\infty}(\Omega) $ and a family $(\underline{u}_{\delta})  \subset W^{1,\infty}_0(\Omega)$ such that:
\begin{equation}\label{weak supersolution}
\int_{\Omega}\nabla \overline{u}\nabla vdx \geq \int_{\Omega}\frac{1}{m}f(x, \overline{u}, \nabla \overline{u})vdx ~~ \forall v\in H_{0}^{1}(\Omega ), \; v\geq 0  ~~~\mbox{and} ~~ \overline{u} \geq 0 ~~ \mbox{on} ~~ \partial \Omega,
\end{equation}
$$
\|\underline{u}_\delta\|_{1,\infty} \to 0 ~~ \mbox{as} ~~ \delta \to 0,
$$
$$
\underline{u}_\delta \leq \overline{u} ~~ \mbox{in} ~~ \Omega ~~ \mbox{for} ~~ \delta ~~ \mbox{small enough},
$$
and given $\alpha >0$, there is $\delta_0>0$ such that

\begin{equation}\label{weak subsolution}
\int_{\Omega}\nabla \underline{u}_\delta \nabla vdx \leq \int_{\Omega}\frac{1}{\alpha}f(x, \underline{u}_{\delta}, \nabla \underline{u}_\delta)vdx ~~ \forall v\in H_{0}^{1}(\Omega ), \; v\geq 0  ~~~\mbox{for} ~~ \delta \leq \delta_0 .
\end{equation}
Then there is a small enough $\delta >0$ such that problem $(P)$ has a weak solution $u$ satisfying $\underline{u}_\delta \leq u \leq  \overline{u}$.
\end{thm}

\section{Preliminary Results}

In this section we introduce some concepts and results in order to attack problem $(P)$. The abstract results concerning monotone operators can be found, for instance, in Lions \cite{lions}, Ne$\check{c}$as \cite{necas} and Pascali \& Sburlan \cite{pascali}

\begin{definition}\label{monotone operator}
Let $E$ be a reflexive Banach space and $E^{\ast}$ its topological dual. A nonlinear mapping $A:D(A)\subset E\rightarrow E^{\ast}$
is said to be monotone if it satisfies
\begin{equation}\label{inequality monotone}
\langle Au-Av,u-v\rangle \geq 0 \;\; u,v\in D(A).
\end{equation}
If the inequality (\ref{inequality monotone}) is strict for $u\neq v$, we say that $A$ is strict monotone. Here, $\langle \cdot , \cdot \rangle$ means the duality pairing between $E^{\ast}$ and $E$.
\end{definition}

\begin{definition}\label{gradient}
If $E$ is a Hilbert space and $\phi : E\rightarrow \Re$ is $C^{1}$-functional, the gradient of $\phi $,  denoted by $\nabla \phi :E\rightarrow E$, is defined, through the Riesz Representation Theorem, by
$$
\langle \nabla \phi (u), w \rangle =\phi' (u)w \;\; \forall u, w \in E,
$$
where $\langle \cdot \; , \; \cdot \rangle $ is the inner product in $E$.
\end{definition}

\begin{lem}\label{convex-monotone}
If $E$ is a Hilbert space and $\phi \in C^{1}(E, \Re )$, then $\phi$ is convex (strictly convex) if, and only if, $\nabla \phi$ is monotone (strictly monotone).
\end{lem}

\begin{definition}\label{S+}
Let $E$ be a Banach space and $\mathcal{C}\subset E$ a closed convex set. An operator $T: \mathcal{C}\rightarrow E^{\ast}$ is said to be of type $(S_{+})$ provided that whenever $x_{n}\rightharpoonup x$ in $E$ and
\begin{equation}\label{condition S+}
\limsup_{n \to +\infty} \; \langle Tx_{n}, x_{n}-x\rangle \leq 0,
\end{equation}
then $x_{n}\rightarrow x$ in $E$.

We remark that the condition (\ref{condition S+}) can be rewritten as
\begin{equation}
\limsup_{n \to +\infty} \; \langle Tx_{n}-Tx, x_{n}-x\rangle \leq 0.
\end{equation}
\end{definition}

\begin{definition}\label{pseudomonotone definition}
Let $E$ be a Banach space and $B: E\rightarrow E^{\ast}$ an operator. We say that $B$ is pseudomonotone if $u_{n}\rightharpoonup u$ in $E$ and
\begin{equation}
\limsup_{n \to +\infty} \; \langle Bu_{n},u_{n}-u \rangle \leq 0,
\end{equation}
then
\begin{equation}
\liminf_{n \to +\infty} \;  \langle Bu_{n},u_{n}-v \rangle \geq \langle B(u),u-v \rangle \;\; \forall v\in E.
\end{equation}
\end{definition}

\begin{definition}\label{demicontinuous definition}
We say that $T:E\rightarrow E^{\ast}$ is demicontinuous if $x_{n}\rightarrow x$ in $E$ implies that $Tx_{n}\rightharpoonup Tx$ in $E^{\ast}$.
\end{definition}

\begin{lem}\label{demicontinuous S+}
Any demicontinuous operator $T:E\rightarrow E^{\ast}$ of type $(S_{+})$ is pseudomonotone.
\end{lem}

\begin{thm}\label{B surjective}
Let $E$ be a reflexive and separable Banach space and $B:E \rightarrow E^{\ast}$ an operator satisfying
\begin{description}
  \item[(i)] $B$ is coercive, i.e.,
  \begin{equation}\label{B coercive}
  \frac{\langle B(u),u\rangle}{\|u\|}\rightarrow +\infty \;\;\mbox{as}\;\; \|u\|\rightarrow +\infty \;
  \end{equation}
  \item[(ii)] $B$ is bounded and continuous;
  \item[(iii)] $B$ is pseudomonotone.
\end{description}
Then $B$ is surjective, that is, $B(E)=E^{\ast}$.
\end{thm}

Next, $\| \cdot \|$ will denote the usual norm $\|u\|= \left(\dis\int_{\Omega}|\nabla u|^{2}dx\right)^{\frac{1}{2}}$ in $H_{0}^{1}(\Omega )$.

\begin{lem}\label{L monotone}
The operator $L:H_{0}^{1}(\Omega )\rightarrow H^{-1}(\Omega )$ given by
\begin{equation}\label{L monotone}
\langle Lu,v\rangle =\dis\int_{\Omega}M(\|u\|^{2})\nabla u \nabla vdx
\end{equation}
is strictly monotone.
\end{lem}

\dem Let us consider $G:H_{0}^{1}(\Omega ) \rightarrow \Re$ given by
\begin{equation}\label{G definition}
G(u)=\frac{1}{2}\widehat{M}(\|u\|^{2})\;\; \forall u \in H_{0}^{1}(\Omega ),
\end{equation}
where $\widehat{M}(t)=\dis\int_{0}^{t}M(\tau )d\tau$. Because $M$ is positive and continuous, we have that $G$ is strictly convex. Furthermore
\begin{equation}\label{G' = L}
G'(u)v=\langle \nabla G(u),v\rangle \dis\int_{\Omega}M(\|u\|^{2})\nabla u \nabla vdx=\langle Lu,v\rangle \;\; \forall u,v\in H_{0}^{1}(\Omega ),
\end{equation}
that is, $\nabla G =L$ and so, in view of Lemma \ref{convex-monotone}, $L$ is strictly monotone. \fim

\begin{lem}\label{L is S+}
$L$ is of type $(S_{+})$.
\end{lem}

\dem Let $(u_{n})$ be a sequence in $H_{0}^{1}(\Omega )$ such that
\begin{equation}\label{weak convergence 1}
u_{n} \rightharpoonup u \;\; \mbox{in} \;\; H_{0}^{1}(\Omega)
\end{equation}
and
\begin{equation}\label{limsup1}
\limsup_{n \to +\infty} \; \langle Lu_{n}, u_{n}-u\rangle \leq 0.
\end{equation}

We have to prove that $u_{n}\rightarrow u$ in $H_{0}^{1}(\Omega )$. For this, we first note that
$$
\langle Lu_n,u_n-u\rangle    =  M(\|u_n\|^{2})\dis\int_{\Omega}|\nabla u_n|^{2}dx-M(\|u_n\|^{2})\dis\int_{\Omega}\nabla u_n\nabla udx
$$
that is,
$$
 \frac{1}{M(\|u_n\|^{2})} \langle Lu_n,u_n-u\rangle  =  \dis\int_{\Omega}|\nabla u_n|^{2}dx-\dis\int_{\Omega}\nabla u_n\nabla udx
$$
Note that $M(\|u_n\|^{2})\geq m>0$, and so,
$$
0\geq \limsup_{n \to \infty} \|u_n\|^{2}-\|u\|^{2},
$$
which implies
$$
\|u\|^{2}\geq \limsup_{n \to \infty} \|u_n\|^{2} \geq \liminf_{n \to \infty} \|u_n\|^{2} \geq \|u\|^{2},
$$
from where it follows that $\|u_n\|^{2}\rightarrow \|u\|^{2}$. Invoking the weak convergence $u_n\rightharpoonup u$ in $H_{0}^{1}(\Omega )$, we see that $u_n\rightarrow u$ in $H_{0}^{1}(\Omega )$, and the proof of the lemma is over. \fim

\section{Proof of the Main Theorem}

From now on, we fix $R>0$ large enough such that
$$
\|\nabla \overline{u}\|_{\infty}, \|\nabla \underline{u}_{\delta}\|_{\infty} \leq R
$$
for all $\delta$ small enough, where $\overline{u}$ and $u_\delta$ were given in Theorem \ref{main result}. We recall that if $\overrightarrow{V}=(V_{1}, \ldots , V_{N})\in (L^{\infty}(\Omega))^{N}$, we have $\|\overrightarrow{V}\|_{\infty}=\dis\max_{1\leq i\leq N}\|V_{i}\|_{\infty}$.  Moreover, we set the function $g_R:\mathbb{R} \to \mathbb{R}$ given by
$$
g_R(t)=
\left\{
\begin{array}{l}
t, ~~ \mbox{if} ~~ |t| \leq R,\\
\mbox{}\\
R, ~~ \mbox{if} ~~  t \geq R, \\
\mbox{}\\
-R, ~~ \mbox{if} ~~ t \leq -R.
\end{array}
\right.
$$
Here, we would like to point out that
\begin{equation} \label{E0}
g_R(t)=t ~~ \mbox{if} ~~ |t| \leq R
\end{equation}
and
$$
|g_{R}(t)|= \min\{R, |t|\} ~~ \mbox{for all} ~~ t \in \mathbb{R}.
$$
Hence,
\begin{equation} \label{E1}
|g_{R}(t)|\leq R ~~~\mbox{and} ~~ |g_{R}(t)|\leq |t| ~~~\mbox{for all} ~~ t \in \mathbb{R}.
\end{equation}

 Taking into account the above function $g_{R}$ and and their properties, we will consider the following auxiliary function $f_R:\Omega \times \mathbb{R} \times \mathbb{R}^{N} \to [0,+\infty)$ given by
$$
f_R(x,t,y)=f(x,t,\stackrel{\rightarrow}{g_R}(y)),
$$
where $\stackrel{\rightarrow}{g_R}(y)=(g_R(y_1),g_R(y_2),...,g_R(y_N))$. Using the definition of the function $f_R$, it follows the ensuing estimates:
\begin{equation} \label{E2}
|f_R(x,t,y)| \leq h(x,t)(1+|\stackrel{\rightarrow}{g_R}(y)|^{\eta}) \leq  h(x,t)(1+R^{\eta}N^{\frac{\eta}{2}})
\end{equation}
and
\begin{equation} \label{E3}
|f_R(x,t,y)| \leq h(x,t)(1+|\stackrel{\rightarrow}{g_R}(y)|^{\eta}) \leq  h(x,t)(1+|y|^{\eta}).
\end{equation}
Furthermore, it is crucial observing that
\begin{equation}  \label{E4}
f_R(x,t,y)=f(x,t,y) ~~~\mbox{if} ~~ |y| \leq R,
\end{equation}
and so,
$$
f_R(x,\overline{u},\nabla \overline{u})= f(x,\overline{u},\nabla \overline{u}) ~~~\mbox{and} ~~ f_R(x,{u}_\delta,\nabla \overline{u}_\delta)= f(x,\overline{u}_\delta,\nabla \overline{u}_\delta).
$$

Using function $f_R$, we are able to fix the following auxiliary problem
$$
\left\{
\begin{array}{rclcc}
  -M\left(\int_{\Omega}|\nabla u|^{2}dx\right)\Delta u & = & f_{R}(x,u, \nabla u ) & \mbox{in} & \Omega , \\

  u  & = & 0 & \mbox{on} & \partial \Omega.
\end{array}
\right.
\eqno{(AP)}
$$

Our intention is proving the existence of a solution $u_R$ for $(AP)$ with $\|\nabla u_R\|_{\infty} \leq R$ if $R$ is large enough and, because of  (\ref{E4}), we can guarantee  that $u_R$ is a solution of the original problem $(P)$.

\subsection{Supersolution}

In this subsection, we will be concerned on supersolutions of the problem $(AP)$.

\begin{definition}\label{supersolution definition}
We say that $w \in W^{1, \infty}(\Omega)$ is a supersolution of the problem $(AP)$ if
\begin{equation}\label{inequality supersolution}
-M(\|w\|^{2})\Delta w\geq f_R(x, w, \nabla w) \;\;\mbox{in}\;\; \Omega ~~ \mbox{and} ~~ w \geq 0 ~~ \mbox{on} ~~ \partial \Omega,
\end{equation}
\end{definition}
in the weak sense, that is,
\begin{equation}\label{weak supersolution}
M(\|w\|^{2})\int_{\Omega}\nabla w\nabla vdx \geq \int_{\Omega}f_R(x, w, \nabla w)vdx
\end{equation}
$\forall v\in H_{0}^{1}(\Omega ), \; v\geq 0$ a.e. in $\Omega$

How to get a supersolution to the problem $(AP)$?  Under the hypotheses of Theorem \ref{main result}, we know that $\overline{u}\in W^{1, \infty}(\Omega)$ verifies
\begin{equation}\label{local supersolution}
-\Delta \overline{u} \geq \frac{1}{m}f(x, \overline{u}, \nabla \overline{u}) \;\;\mbox{in}\;\; \Omega.
\end{equation}

Since $M(t)\geq m>0$ for all $t\geq 0$ and $f(x, \overline{u}, \nabla \overline{u})=f_R(x, \overline{u}, \nabla \overline{u})$, we deduce that
$\overline{u}\in W^{1, \infty}(\Omega)$ is a supersolution of the problem $(AP)$.

We point out that sub and supersolutions for quasilinear  local problems like
\begin{equation} \label{local problem}
\left\{
\begin{array}{rclcc}
  -\Delta u & = & f(x,u, \nabla u ) & \mbox{in} & \Omega , \\

  u  & = & 0 & \mbox{on} & \partial \Omega
\end{array}
\right.
\end{equation}
were studied in \cite{cuestaleon}.

\begin{lem}\label{a priori supersolution}
Let $u_R \in H_{0}^{1}(\Omega )\cap L^{\infty}(\Omega )$ be a weak  solution of the  problem $(AP)$ with $0<u_R \leq \overline{u}$ a.e. in $\Omega$. Then there is a constant $K=K(\|\overline{u}\|_{\infty},R)$ such that
\begin{equation}
\|u_R\|^{2}\leq K.
\end{equation}
\end{lem}

\dem Setting $T= \|\overline{u}\|_{\infty}$, by condition $(f_{1})$ combined with (\ref{E2}), there is a constant  $C=C(T)>0$ such that
\begin{equation}\label{consequence f1}
|f_R(x,t,y)|\leq C(1+R^{\eta}N^{\frac{\eta}{2}})=C_1
\end{equation}
for all $(x,t,y)\in \Omega \times [0,T]\times \Re^{N}$. Since $u_R$ is a solution of $(AP)$, we have
\begin{equation}
M(\|u_R\|^{2})\dis\int_{\Omega}\nabla u_R \nabla vdx=\dis\int_{\Omega}f_R(x,u_R, \nabla u_R)vdx \;\; \forall v\in H_{0}^{1}(\Omega )
\end{equation}
and so
\begin{equation}
M(\|u_R\|^{2})\|u_R\|^{2}=\dis\int_{\Omega}f_R(x,u_R, \nabla u_R)u_Rdx.
\end{equation}
Invoking (\ref{consequence f1}), we obtain
\begin{equation}
M(\|u_R\|^{2})\|u_R\|^{2}\leq C_1\dis\int_{\Omega}|u_R|dx \leq C_1\dis\int_{\Omega}|\overline{u}|dx
\end{equation}
leading to
\begin{equation}
m\|u_R\|^{2}\leq C_2,
\end{equation}
from where it follows that there is $K>0$ satisfying $\|u_R\|^{2}\leq K$. \fim

\subsection{Subsolution}\label{subsolution section}

In this section we will be concerned on subsolutions of $(AP)$.

\begin{definition}\label{subsolution definition}
We say that $w \in W^{1, \infty}(\Omega)$ is a subsolution of the problem $(AP)$ if
\begin{equation}\label{inequality subsolution}
-M(\|w\|^{2})\Delta w \leq f_R(x, w, \nabla w) \;\;\mbox{in}\;\; \Omega ~~ \mbox{and} ~~ w \leq 0 ~~~\mbox{on} ~~ \partial \Omega,
\end{equation}
\end{definition}
in the weak sense, that is,
\begin{equation}\label{weak subsolution}
M(\|w\|^{2})\int_{\Omega}\nabla w\nabla vdx \leq \int_{\Omega}f_R(x, w, \nabla w)vdx
\end{equation}
$\forall v\in H_{0}^{1}(\Omega ), \; v\geq 0$ a.e. in $\Omega$

In order to construct a subsolution,  we consider the family $(\underline{u}_{\delta}) \subset W_{0}^{1, \infty}(\Omega )$ mentioned in Theorem \ref{main result}, we know that there is $\delta^{*}>0$ such that
 \begin{equation}
\underline{u}_{\delta}\leq \overline{u} \;\; \forall \delta \in [0, \delta^{\ast}],
\end{equation}
with
\begin{equation}
\|\underline{u}_{\delta}\|_{1, \infty}\rightarrow 0 \; \mbox{as}\; \delta \rightarrow 0^{+}.
\end{equation}

Thereby, fixing $\alpha =\displaystyle\max_{t\in [0, 1 ]}M(t)$, we can reduce if necessary $\delta^{*}$ to get
\begin{equation}
-\Delta \underline{u}_{\delta}\leq \frac{1}{\alpha}f(x, \underline{u}_{\delta}, \nabla \underline{u}_{\delta}) ~~ \mbox{in} ~~ \Omega ~~ \mbox{and} ~~ u_\delta=0 ~~ \mbox{on} ~~ \partial \Omega.
\end{equation}
Once that  $f(x, \underline{u}_{\delta}, \nabla \underline{u}_{\delta})=f_R(x, \underline{u}_{\delta}, \nabla \underline{u}_{\delta})$, we can claim that  $\underline{u}=u_\delta$ for $\delta \in (0, \delta^{*})$ is a subsolution of $(AP)$.

%%%%%%%%%%%%%%%%%%%%%%%%%%%%%%%%%%%%%%%%%%%%%%%%%%%%%%%%%%%%%%%%%%%%%%%%%%%%%%%%%%%%%%%%%%%%%%%%%%%%%%%%%
%%%%%%%%%%%%%%%%%%%%%%%%%%%%%%%%%%%%%%%%%%%%%%%%%%%%%%%%%%%%%%%%%%%%%%%%%%%%%%%%%%%%%%%%%%%%%%%%%%%%%%%%%%
\subsection{Another Auxiliary Problem}

In what follows, we define
$$
z_R(x,t,y)=\left\{
\begin{array}{rl}
  f_R(x,\underline{u}(x),\nabla \underline{u}(x)), & t\leq \underline{u}(x), \\
  f_R(x,t,y), & \underline{u}(x)\leq t \leq \overline{u}(x), \\
  f_R(x, \overline{u}(x), \nabla \overline{u}(x)), & t\geq \overline{u}(x)
\end{array}
\right.
$$
and for $l\in (0,1)$ we define the function
$$
\gamma_R (x,t)=-(\underline{u}(x)-t)^{l}_{+}+(t-\overline{u}(x))_{+}^{l}.
$$

Using the above functions, we consider below a second auxiliary problem
\begin{equation} \label{auxiliary problem}
\left\{
\begin{array}{rclcc}
  -M\left(\dis\int_{\Omega}|\nabla u|^{2}dx\right)\Delta u & = & z_R(x,u, \nabla u )-\gamma_R (x,u) & \mbox{in} & \Omega , \\

  u  & = & 0 & \mbox{on} & \partial \Omega.
\end{array}
\right.
\end{equation}

Next, our goal is  proving the existence of a solution for the problem (\ref{auxiliary problem}). To this end, we will use Theorem \ref{B coercive} to the operator
$$
\begin{array}{rcl}
  B: H_{0}^{1}(\Omega ) & \rightarrow & H^{-1}(\Omega ) \\
  u & \mapsto & B(u)
\end{array}
$$
where
$$
\begin{array}{rcl}
  B(u): H_{0}^{1}(\Omega ) & \rightarrow & \Re \\
  v & \mapsto & \langle B(u),v\rangle
\end{array}
$$
is given by
$$
\langle B(u),v\rangle =M(\|u\|^{2})\int_{\Omega}\nabla u \nabla v dx -\int_{\Omega}z_R(x,u, \nabla u)vdx +\int_{\Omega}\gamma_R (x,u)vdx.
$$
 In what follows, we are going to show  that $B$ is onto. So, there exists $u_R\in H_{0}^{1}(\Omega )$ such that $B(u_R)=0$ in $H^{-1}(\Omega )$. Consequently, $u_R$ is a weak solution of the auxiliary problem. If such a solution enjoys $\underline{u}\leq u_R \leq \overline{u}$ a.e. in $\Omega$ we get a solution of problem $(AP)$.

Plainly  $B$ is continuous. In what follows, we fix our attention to others properties of $B$ in order to apply  Theorem \ref{B coercive}.

\begin{lem}
$B$ is coercive.
\end{lem}

\dem First note that
$$
\langle B(u),u\rangle =M(\|u\|^{2})\|u\|^{2}-\int_{\Omega }z_R(x,u, \nabla u )udx +\int_{\Omega }\gamma_R (x,u)udx.
$$
It follows from the definition of $z_R$ that there exists $C=C(R)>0$ such that
$$
z_R(x,t,y)\leq C ~~~ \forall (x,t,y) \in \Omega \times \mathbb{R} \times \mathbb{R}^{N}
$$
and
$$
|\gamma (x,t)|\leq C_{1}+C_{2}t^{l} ~~~ \forall (x,t) \in \Omega \times \mathbb{R}.
$$
Consequently,
$$
|z_R(x,u, \nabla u)u| \leq C|u|
$$
and
$$
|\gamma (x,u)|\leq C_{1}|u|+C_{2}|u|^{l+1}.
$$
From these last inequalities,
$$
-\int_{\Omega}z_R(x,u, \nabla u)udx \geq -\int_{\Omega}|z_R(x,u, \nabla u)|udx \geq -C_{1}\|u\|
$$
and
$$
\int_{\Omega}\gamma (x,u)udx \geq -\int_{\Omega}|\gamma (x,u)u|dx \geq -C_{3}\int_{\Omega}|u|dx -C_{4}\int_{\Omega}|u|^{l+1}dx,
$$
that is,
$$
\int_{\Omega}\gamma (x,u)udx \geq -C_{5}\|u\|-C_{6}\|u\|^{l+1}.
$$

Since $M(t)\geq m>0$ for all $t\geq 0$, one has
$$
\langle B(u),u \rangle \geq m \|u\|^{2}-C_{7}\|u\|-C_{6}\|u\|^{l+1}
$$
which yields
$$
\frac{\langle B(u),u \rangle}{\|u\|}\geq m\|u\| -C_{7}-C_{6}\|u\|^{l}
$$
and the result follows because $l\in (0,1)$. \fim

\begin{lem}
$B$ is pseudomonotone.
\end{lem}

\dem Let $(u_{n}) \subset H_{0}^{1}(\Omega )$ be a sequence satisfying
$$
u_{n}\rightharpoonup u ~~ \mbox{in} ~~ H_{0}^{1}(\Omega ) ~~ \mbox{and} ~~  \limsup_{n \to \infty} \; \langle B(u_{n}), u_{n}-u \rangle \leq 0,
$$
and recall that
\begin{equation}
\langle B(u_{n}), u_{n}-u \rangle =\langle L u_{n}, u_{n}-u\rangle -\int_{\Omega}h(x,u_{n}, \nabla u_{n})(u_{n}-u)dx+\int_{\Omega}\gamma (x,u_{n})(u_{n}-u)dx.
\end{equation}
Note that
$$
\left|\int_{\Omega}z_R(x,u_{n}, \nabla u_{n})(u_{n}-u)dx\right|\leq C|u_{n}-u|_{1} \rightarrow 0
$$
and
$$
\int_{\Omega}\gamma_R (x,u_{n})|u_{n}-u|dx \rightarrow 0.
$$
Consequently,
$$
\limsup_{n \to \infty} \; \langle B(u_{n}),u_{n}-u\rangle =\limsup_{n \to \infty} ~\langle Lu_{n}, u_{n}-u\rangle .
$$
Since $L$ is an operator of the type $(S_{+})$, it follows that $u_{n}\rightarrow  u$ in $H_{0}^{1}(\Omega )$ and invoking the continuity of $B$, we obtain
$$
\liminf_{n \to \infty} \; \langle B(u_{n}),u_{n}-u \rangle = \langle B(u),u-v \rangle \;\; \forall v \in H_{0}^{1}(\Omega ),
$$
showing that $B$ is pseudomonotone. \fim

\vspace{0.5 cm}

From the above lemmas, the operator $B$ enjoys all  the hypotheses of Theorem \ref{B coercive} and so $B$ is onto. Consequently, there is $u_R \in H_{0}^{1}(\Omega )$ such that $B(u_R)=0$.

\subsection{Existence of Solution for $(AP)$}

As we remarked before, it is enough to show that $\underline{u} \leq u_R \leq \overline{u}$. In this section, we will denote $u_R$ by $u$.

%First of all, we consider $\underline{u}:=\underline{u}_{\delta}$ with $\delta >0$ small enough.

\vspace{0.5cm}

\noindent ${ \bf 1^{st}}$ {\bf Step}. $u \leq \overline{u}$.

\vspace{0.5cm}

For this first step, we take $v=(u-\overline{u})_{+}$ as a test function. Then,
$$
M(\|u\|^{2})\dis\int_{\Omega}\nabla u \nabla (u-\overline{u})_{+}dx = \dis \int_{\Omega}z_R(x,u, \nabla u)(u-\overline{u})_{+}- \dis \int_{\Omega}\gamma_R (x,u)(u-\overline{u})_{+}
$$
Thus,
$$
\begin{array}{rcl}
  \int_{\Omega}\nabla u \nabla (u-\overline{u})_{+}dx & = & \int_{\Omega}\frac{1}{M(\|u\|^{2})}f_R(x, \overline{u}, \nabla \overline{u})(u-\overline{u})_{+}dx-\frac{1}{M(\|u\|^{2})} \int_{\Omega}(u-\overline{u})_{+}^{l+1}dx \\
	\mbox{}\\
   & \leq  & \frac{1}{m}\int_{\Omega}f_R(x, \overline{u}, \nabla \overline{u})(u-\overline{u})_{+}dx-\frac{1}{M(\|u\|^{2})} \int_{\Omega}(u-\overline{u})_{+}^{l+1}dx \\
	\mbox{}\\
   & \leq  & \int_{\Omega}\nabla \overline{u}\nabla (u-\overline{u})_{+}dx-\frac{1}{M(\|u\|^{2})} \int_{\Omega}(u-\overline{u})_{+}^{l+1}dx.
\end{array}
$$
Combining these inequalities, we get
$$
0\leq \int_{\Omega}|\nabla (u-\overline{u})_{+}|^{2}dx\leq -\frac{1}{M(\|u\|^{2})} \int_{\Omega}(u-\overline{u})_{+}^{l+1}dx \leq 0,
$$
from where it follows that $ u \leq \overline{u}$ in $\Omega$.

\vspace{0.5cm}

\noindent  ${\bf 2^{nd}}$ {\bf Step}. $\underline{u}\leq u$.

\vspace{0.5cm}

Firstly, we point out   that if $\delta >0$ is small enough, there is $\beta^{*}>0$, independent of $\delta$, such that $\|u\|^{2}\leq \beta^{*}$. Indeed, note that
$$
M(\|u\|^{2})\int_{\Omega}|\nabla u|^{2}dx=\int_{\Omega}z_R(x,u,\nabla u)udx-\int_{\Omega}\gamma_R (x,u)udx.
$$
By the first step, $\gamma_R (x,u)=-(\underline{u}-u)_{+}^{l}$. Then,
$$
m\|u\|^{2}\leq C\int_{\Omega}|u|dx+ \int_{\Omega}(\underline{u}-u)_{+}^{l}|u|
$$
This last inequality gives
$$
m\|u\|^{2} \leq C\|u\|+ C\|\underline{u}\|_{\infty}^{l}\|u\|+C\|u\|^{l+1}.
$$
Thereby, there is $\beta^{*}=\beta^{*}(R,m,l)>0$, independent of $\delta>0$ small enough, such that
$$
\|u\|^{2} \leq \beta^*.
$$

In what follows, we reduce $\delta>0$ if necessary, to get
$$
-\Delta \underline{u}\leq \frac{1}{\alpha^{*}} f_R(x,\underline{u}, \nabla \underline{u})
$$
where $\alpha^{*} =\dis \max_{0\leq t\leq \beta^*}M(t)$. Choosing $v=(\underline{u}-u)_{+}$, we obtain
$$
\begin{array}{rcl}
  M(\|u\|^{2})\int_{\Omega}\nabla u \nabla (\underline{u}-u)_{+}dx & = & \int_{\Omega}z_R(x,u,\nabla u)(\underline{u}-u)_{+}dx-\int_{\Omega}\gamma_R (x,u)(\underline{u}-u)_{+}dx \\
   & = & \int_{\Omega}z_R(x,\underline{u}, \nabla \underline{u})(\underline{u}-u)_{+}dx+\int_{\Omega}(\underline{u}-u)_{+}^{l+1}dx
\end{array}
$$
and so
$$
\int_{\Omega}\nabla u\nabla (\underline{u}-u)_{+}dx=\int_{\Omega}\frac{1}{M(\|u\|^{2})}z_R(x,\underline{u}, \nabla \underline{u})(\underline{u}-u)_{+}dx+\int_{\Omega}\frac{1}{M(\|u\|^{2})}(\underline{u}-u)_{+}^{l+1}dx.
$$
Hence
$$
\begin{array}{rcl}
  \int_{\Omega}\nabla u \nabla (\underline{u}-u)_{+}dx & \geq  & \int_{\Omega}\frac{1}{\alpha^{*}}f_R(x, \underline{u} , \nabla \underline{u})(\underline{u}-u)_{+}dx +\frac{1}{M(\|u\|^{2})}\int_{\Omega}(\underline{u}-u)_{+}^{l+1}dx \\
	\mbox{}\\
   & \geq  & \int_{\Omega}\nabla \underline{u}\nabla (\underline{u}-u)_{+}dx+\frac{1}{M(\|u\|^{2})}\int_{\Omega}(\underline{u}-u)_{+}^{l+1}dx.
\end{array}
$$
Then
$$
0\geq \int_{\Omega}|\nabla (\underline{u}-u)_{+}|^{2}dx+\frac{1}{M(\|u\|^{2})}\int_{\Omega}(\underline{u}-u)_{+}^{l+1}dx \geq 0
$$
and this implies that $(\underline{u}-u)_{+}=0$. Thus, $\underline{u}\leq u$ in $\Omega$, and the proof of the existence of solution for $(AP)$ is over.

\subsection{Existence of Solution for $(P)$}

To begin with, we observe that in the last subsection we proved the existence of a solution $u_R$ of $(AP)$ verifying $\underline{u} \leq u_R \leq \overline{u}$ in $\Omega$. Here,  we would like point out that $\underline{u}$ and $\overline{u}$ does not depend of $R$, for $R$ large enough. In what follows, we denote $u_R$ by $u$.

Our goal is to show that there is $R^{*}>0$ such that
$$
\|\nabla u\|_{\infty} \leq R ~~~~ \mbox{for} ~~~ R \geq R^{*}.
$$

 By Elliptic Regularity,
\begin{equation*}
u \in W^{2,p}(\Omega )\,\,\,\forall p\in \lbrack 1,+\infty ),
\end{equation*}
because $f_R \in L^{\infty }([0,+\infty ))$ and $u\in L^{\infty
}(\Omega )$.  From now on, we will
fix $p$ such that
\begin{equation}
W^{2,p}(\Omega )\hookrightarrow C^{1,\alpha }(\overline{\Omega })
\label{CC1}
\end{equation}%
is a continuous embedding. Now, we observe that $u$ is a solution of the
problem
\begin{equation*}
-\Delta u+u=B_{R}(x)(1+|\nabla u|^{2}),
\end{equation*}%
where
\begin{equation*}
B_{R}(x)=\frac{u+\frac{f_R(x,u(x),\nabla u(x))}{M(\|u\|^{2})}}{1+|\nabla
u|^{2} }.
\end{equation*}%
Once that
\begin{equation*}
|f_{R}(x,t,y)|\leq h(x,t)(1+|y|^{\eta}) ~~~\forall (x,t,y) \in \Omega \times \mathbb{R} \times \mathbb{R}^{N},
\end{equation*}%
combining the  fact that $\eta \in [0,2]$,  $\underline{u} \leq u \leq \overline{u}$ in $\Omega$ and $\|\underline{u}\|_\infty, \|\overline{u}\|_\infty$ does not depend  of $R$, for $R$ large enough, the conditions $(f_1)$ and $(M_2)$ guarantee the existence of $C^{\ast }>0$, independent of $R$, such that
\begin{equation*}
|B_{R}(x)|\leq C^{\ast }\,\,\,\forall x\in \Omega, \,\,\,\mbox{for} ~~ R ~~ \mbox{large enough}.
\end{equation*}%
Thereby, there is $R_1>0$ such that
\begin{equation}
\Vert B_{R}\Vert _{\infty }\leq C^{\ast }\,\,\,\forall R>R_1.  \label{CC2}
\end{equation}%
By using a result due to Amann \& Crandall \cite[Lemma 4]{Amann}, there is
an increasing function $\gamma _{0}:[0,+\infty )\rightarrow \lbrack 0,\infty
)$, depending only of $\Omega $, $p$ and $N$, and satisfying
\begin{equation*}
\Vert u\Vert _{W^{2,p}(\Omega )}\leq \gamma _{0}(\Vert B_{R}\Vert _{\infty
}).
\end{equation*}%
Combining the last inequality with (\ref{CC1}) and (\ref{CC2}), we get
\begin{equation*}
\Vert u\Vert _{C^{1,\alpha }(\overline{\Omega} )}\leq C\gamma _{0}(C^{\ast
}),
\end{equation*}%
for some $C>0$. Fixing
\begin{equation*}
K_1=C\gamma _{0}(C^{\ast }),
\end{equation*}%
we derive that
\begin{equation*}
|\frac{\partial u(x)}{\partial x_{i}}|\leq K_1\,\,\,\forall x\in \overline{%
\Omega }\,\,\ \mbox{and}\,\,\,i=1,2,..,N.
\end{equation*}%
Thereby,
\begin{equation*}
|\nabla u(x)|\leq NK_1\,\,\,\forall x\in \overline{\Omega },
\end{equation*}%
implying that
\begin{equation*}
\max_{x\in \overline{\Omega }}|\nabla u(x)|\leq NK_1.
\end{equation*}%
Fixing $R_{2}=NK_1$ and $R \geq R^{*}=\max\{R_1,R_2\}$, it follows that
\begin{equation*}
\max_{x\in \overline{\Omega }}|\nabla u(x)|\leq R,
\end{equation*}%
showing that $u$ is a solution of $(P)$ if $R \geq R^{*}$.

\section{Applications}

In this section, we will present two situations in which our main theorem works.

\vspace{0.5 cm}

\noindent {\bf Application 1:}   Our first application is the following problem

\begin{equation} \label{example1}
\left\{
\begin{array}{l}
  -M\left(\int_{\Omega}|\nabla u|^{2}dx\right)\Delta u = \lambda |u|^{q} + |u|^{p}+\mu |\nabla u|^{q} ~~ \mbox{in} ~~ \Omega , \\
u(x)>0  ~~ \mbox{in} ~~ \Omega \\
  u   =  0 ~~ \mbox{on} ~~  \partial \Omega .
\end{array}
\right.
\end{equation}
where $\lambda$ is a positive parameter, $0<q<1<p<+\infty$  and $M$ verifies conditions $(M_1)-(M_2)$.

Here, we must observe that the above problem is a nonlocal version of a well known result due to  Ambrosetti, Brezis \& Cerami \cite{abc} with an additional gradient term $|\nabla u|^{q}$.

We begin observing that it is easy to find a positive function $\overline{u}$ verifying the inequality
$$
-\Delta \overline{u} \geq \frac{1}{m}(\lambda \overline{u}^{q} + \overline{u}^{p}+\mu |\nabla \overline{u}|^{q})
$$
if $\lambda , \mu$ are small enough. It is enough to follow the ideas found in Ambrosetti, Brezis \& Cerami \cite{abc}. Indeed, let $0<e$ in $\Omega$, $e\in C^{1}(\overline{\Omega})$ the only solution of

\begin{equation} \label{e 1 solution}
\left\{
\begin{array}{rclcc}
  -\Delta e & = & 1 & \mbox{in} & \Omega , \\

  e  & = & 0 & \mbox{on} & \partial \Omega .
\end{array}
\right.
\end{equation}

We now take $S>0$ such that
\begin{equation}\label{S large}
m\geq \frac{1}{S^{1-q}}(\lambda \|e\|_{\infty}^{q}+\mu \||\nabla e|\|_{\infty}^{q})+S^{p-1}\|e\|_{\infty}^{p}.
\end{equation}
A straightforward computation shows that there is $0< \lambda^{\ast}$ such that for $0<\lambda , \mu <\lambda^{\ast}$ there is $S>0$ such that the inequality (\ref{S large}) holds true. Hence we can take $\overline{u}:= Se\in W^{1, \infty}(\Omega )$, $S$ as above, so that the first inequality in the Theorem \ref{main result} is satisfied.

Now, fixed $\lambda , \mu  >0$ as before, we consider the family $(u_\delta)$ with $u_\delta=\delta \varphi_1$, $\varphi_1$ is a positive  eigenfunction associated with the principal eigenvalue $\lambda_1$ of $(-\Delta , H_{0}^{1}(\Omega))$.  A simple computation also gives
for all $\alpha >0$ fixed, there exist $\delta^*>0$ such that
$$
-\Delta u_\delta \leq \frac{1}{\alpha}(\lambda u_\delta^{q}+u_\delta^{p}+\mu |\nabla u_{\delta}|^{q}) ~~~\mbox{in} ~~ \Omega .
$$
As it is well known, we can consider $\delta >0$ sufficiently  small such that $u_{\delta}\leq \overline{u}$.

From the above commentaries, we can apply Theorem \ref{main result} to prove the existence of a weak solution $u$ for (\ref{example1}) satisfying $u_{\delta}\leq u\leq \overline{u}$.

\vspace{0.5 cm}

\noindent {\bf Application 2:} ~~ Our next application is concerning the problem

\begin{equation} \label{example2}
\left\{
\begin{array}{l}
  -M\left(\int_{\Omega}|\nabla u|^{2}dx\right)\Delta u  =  Au^{q}(B-u)+|\nabla u|^{\eta} ~~  \mbox{in} ~~  \Omega , \\

  u   =  0 ~~ \mbox{on} ~~\partial \Omega
\end{array}
\right.
\end{equation}
where $A,B$ are positive constants satisfying some properties which will be established later, $\eta  \in (1,2]$ and $q \in (0,1)$. We will find a solution $u$ satisfying $0<u\leq B$ in $\Omega$. First of all, let us consider the continuous function $f:\Omega \times \Re \times \Re^{N}\rightarrow [0,+\infty)$ defined as
$$
f(x,t,y)=
\left\{
\begin{array}{l}
|y|^{\eta}, ~~ \mbox{if} ~~ t \geq B \\
\mbox{}\\
 At^{q}(B-t)+|y|^{\eta}, ~~ \mbox{if} ~~0 \leq  t \leq B \\
\mbox{}\\
|y|^{\eta} , ~~ \mbox{if} ~~ t \leq 0.
\end{array}
\right.
$$
It is clear that the function $\overline{u}\equiv B$ belongs to $W^{1, \infty}(\Omega )$ and satisfies the assumption (\ref{weak supersolution}) in the Theorem \ref{main result}.

If $\lambda_{1}$ is the principal eigenvalue of $(-\Delta , H_{0}^{1}(\Omega ))$ associated to the eigenfunction $\varphi_{1}>0$ in $\Omega$, for each $\alpha >0 $,  there is $\delta^* >0$  such that
$$
\lambda_{1}\delta \varphi_{1}\leq \frac{A}{\alpha}(\delta \varphi_{1})^q(B-\delta\varphi_{1})+\frac{1}{\alpha}|\nabla (\delta\varphi_{1})|^{\eta} ~~~~ \forall \delta \in (0, \delta^*].
$$
Taking $\underline{u}_{\delta}:= \delta \varphi_{1}$ we get $\|\underline{u}_{\delta} \|_{1, \infty}\rightarrow 0$ as $\delta \rightarrow 0$, $\underline{u}_{\delta} \leq \overline{u}\equiv B $ in $\Omega$, if $\delta >0$ is small enough and a straightforward calculation shows that the inequality (\ref{weak subsolution}) holds true. Hence, for $\delta$ sufficiently small, problem (P) possesses a weak solution $u$ satisfying $\underline{u}_{\delta}\leq u \leq B$. Consequently, such a function is a solution of the problem (\ref{example2}).

\begin{rmk}
For some applications concerning the quasilinear problem (P), with $M\equiv 1$, still using a sub and supersolution approach, the reader may consult    Xavier \cite{xavier} and the references therein.
\end{rmk}

\end{document}